\documentclass[preprint]{elsarticle}
\usepackage{amsthm,amsmath}
\usepackage{amsfonts, graphicx,url}

\begin{document}

\journal{Computers and Mathematics with Applications}

\begin{frontmatter}
\title{Conjecture concerning a completely monotonic function}
\author[RISC]{E. Shemyakova}
\author[Khas]{S.I. Khashin}
\author[UWO]{D.J. Jeffrey}
\address[RISC]{Research Institute for Symbolic Computation, Johannes Kepler University,
 Linz, Austria}
\address[Khas]{Mathematics Department, Ivanovo State University, Ivanovo, Russia}
\address[UWO]{Department of Applied Mathematics, The University of Western Ontario,\\
London, Ontario, Canada}

\begin{abstract}
Based on a sequence of numerical computations, a conjecture is presented
regarding the class of functions $H(x;a)=\exp(a)-(1+a/x)^x$, and the open problem of
determining the values of $a$ for which the functions are completely monotonic with respect to $x$.
The critical value of $a$ is
determined here to sufficient accuracy to show that it is not a simple symbolic quantity.
\end{abstract}

\end{frontmatter}

\section{Introduction}

A completely monotonic (CM) function is an infinitely differentiable function $f:(0,\infty)\to \mathbb{R}$ whose derivatives satisfy
\[   \forall x>0,\forall n\in \mathbb{N}, (-1)^n f^{(n)}(x) \ge 0 \ .
\]
The definition is due to Hausdorff~\cite{Hausdorff:1921}, and a collection of important properties can be found in \cite{Widder:1941}.
Alzer and Berg~\cite{Alzer:2002} state the following open problem.
Consider the function
\begin{equation}\label{eq:def}
   H(x;a) = e^a - \left(1+\frac{a}{x}\right)^{x}\ .
\end{equation}
For what values of the parameter $a$ is $H(x;a)$ completely monotonic (CM)?

A problem that is superficially similar to this problem has a remarkably simple solution~\cite{Alzer:2002}.
The function
\begin{equation}\label{eq:known}
   J(x;a,b) = \left(1+\frac{a}{x}\right)^{x+b}-e^a
\end{equation}
is CM if and only if $a\le 2b$.
Although $H(x;a)=-J(x;a,0)$, the problems for $J$ and $H$ are essentially different, because a completely monotonic function must be positive and decreasing.
A more completely descriptive name might be `completely monotonically decreasing function',
but the shorter name is now standard.
In spite of the beautiful result for $J$, a consequence of the conjecture advanced here is that
no similarly simple result applies to $H(x;a)$.

\newcommand{\ac}{a_\mathrm{c}}
\newcommand{\af}{a_\mathrm{f}}
\newcommand{\xf}{x_\mathrm{f}}

Alzer and Berg~\cite{Alzer:2002} show that $H(x;1)$ is CM.
Moreover, Berg~\cite{Berg:2005} has shown, through an equivalent problem, that $H(x;3)$ is not CM.
Therefore, it has been conjectured
that there exists a value $\ac$ such that $H(x;a)$ is completely monotonic for
all $a<\ac$ and not for $a>\ac$.
We obtain here experimentally an estimate for its value, namely $\ac\approx 2.299656443$.

In the field of experimental mathematics, there are several tools for identifying a
real number with an expression containing known mathematical constants. We have applied these
tools to $\ac$, but without success.

\section{Heuristic description of method}
In this section we give a graphical and heuristic description of the method we use.
For convenience, we define the function
\begin{equation}
f(x,a,n) = (-1)^n \frac{\partial^n}{\partial x^n} H(x; a) \ ,
\end{equation}
and consider its properties for various fixed $n$.

\begin{figure}
\begin{center}
(a)\raisebox{-4cm}{
\scalebox{0.27}{\includegraphics{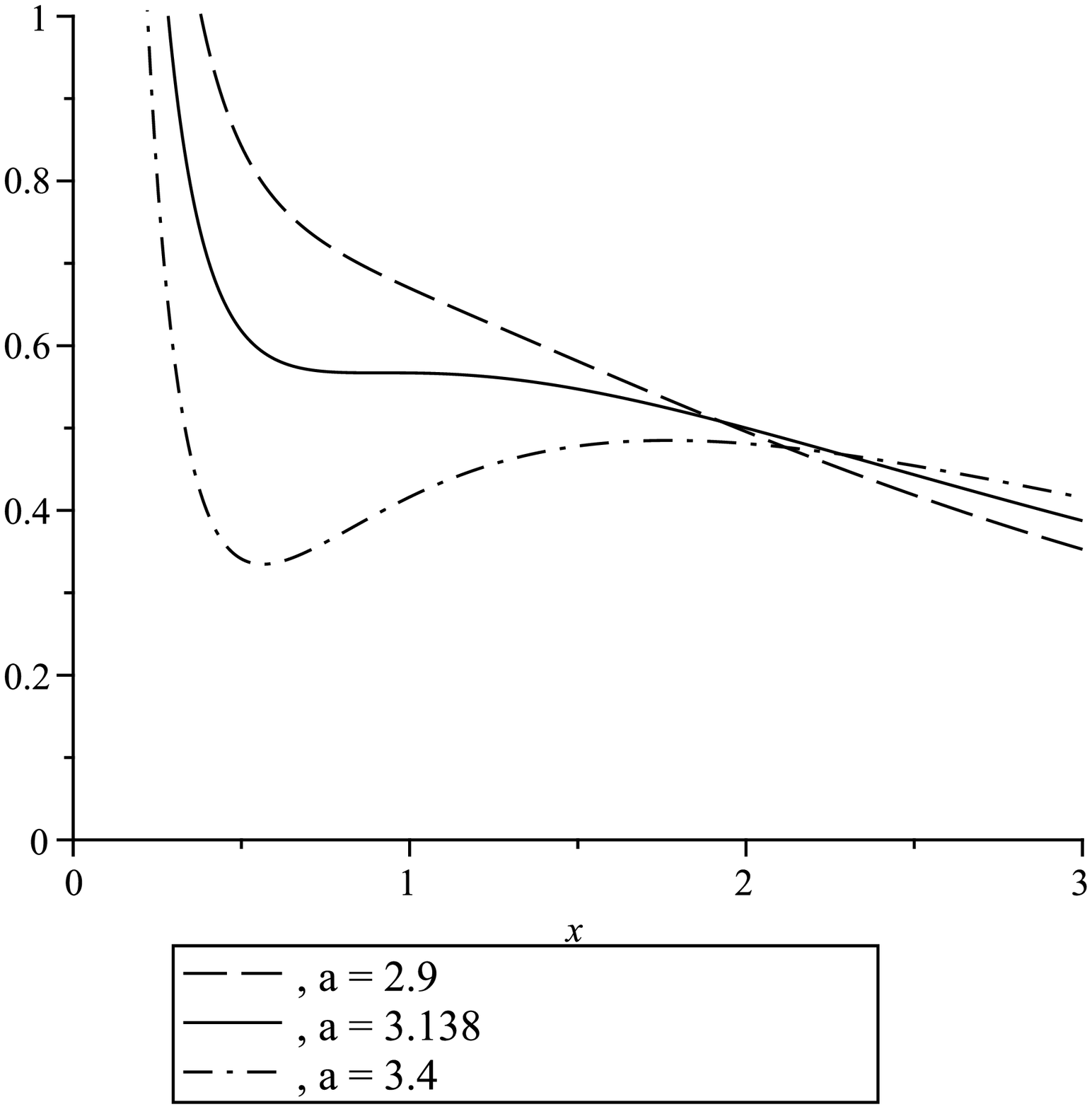}}
}
(b)\raisebox{-4cm}{
\scalebox{0.27}{\includegraphics{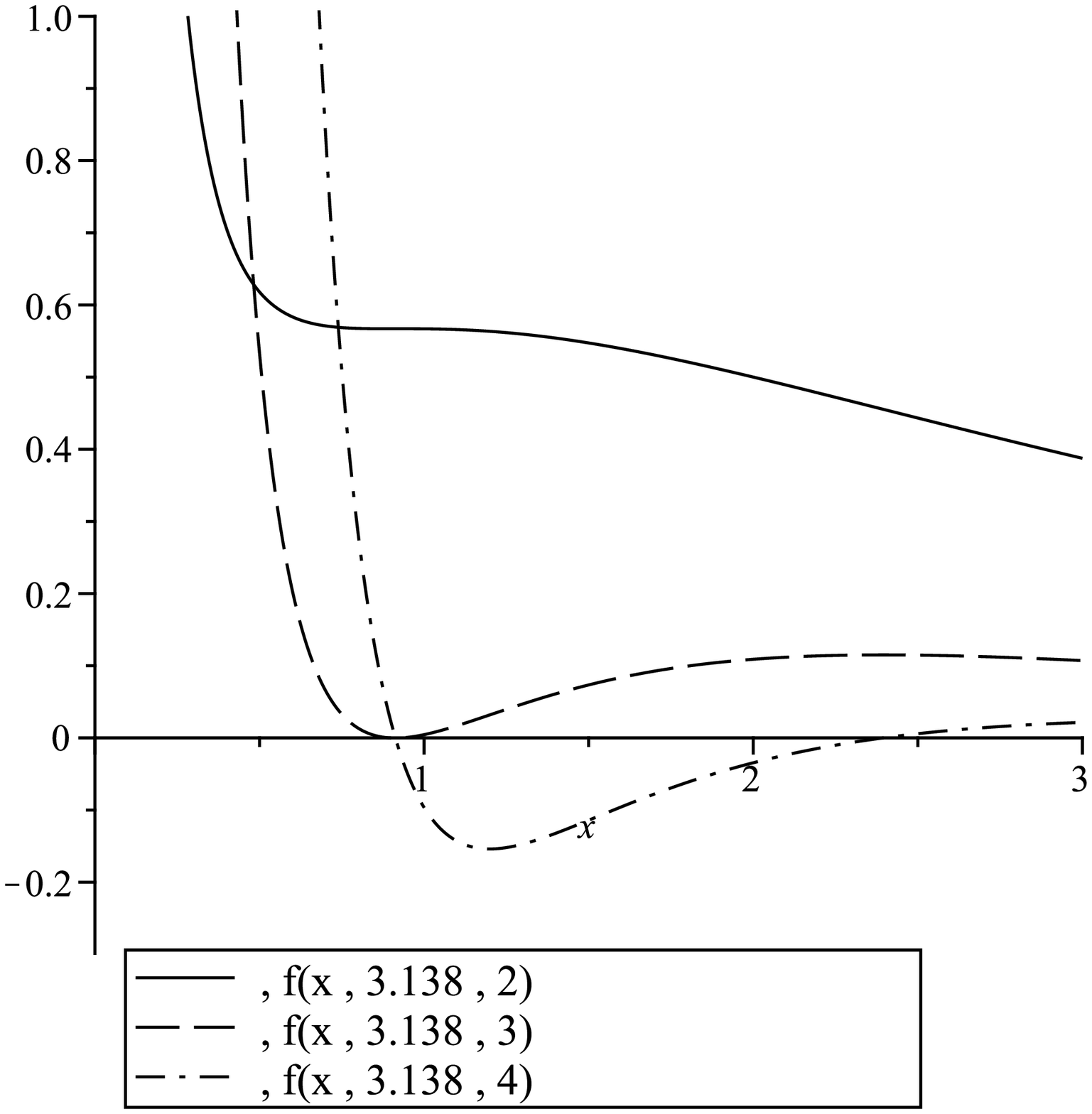}}
}
\end{center}
\caption{\label{fig:1} (a) $f(x,a,2)$ for three values of $a$.
The function $f(x,2.9,2)$ is monotonic, $f(x,3.4,2)$ is not monotonic, and $f(x,3.138,2)$
is the transitional case for $n=2$.
(b) The right plot shows for the critical value $a=3.138$ the successive derivatives $f(x,3.138,2)$, $f(x,3.138,3)$ and $f(x,3.138,4)$.
 }
\end{figure}

We start with the case $n=2$.
Figure \ref{fig:1}(a) shows plots of $f(x,a,2)$, or equivalently the second derivative of $H$.
For $a=2.9$, the function is clearly monotonic decreasing; for $a=3.4$, the function is
clearly not monotonic.
The transition occurs for $a=\af \approx 3.138$, when $f(x,\af,2)$ has an inflexion point at $x=\xf \approx 0.913$, as shown.

The transitional case is further considered in figure \ref{fig:1}(b), where the second, third and fourth derivatives are plotted. As could be anticipated from the shape of
$f(x,\af,2)$, its derivative $f(x,\af,3)$ touches the $x$-axis at $x=\xf$.
Further, $f(x,\af,4)$ crosses the axis at the same point.
From these observations, it is obvious that the pair of values $\af,\xf$ satisfy the equations
\begin{eqnarray}
  f(\xf,\af,3) &=& 0\ , \\
  f(\xf,\af,4) &=& 0\ .
\end{eqnarray}

Now consider the case $n=3$. The above series of observations can be repeated for $f(x,a,3)$. Again there is a range of values of $a$ for which the function is monotonic decreasing, and there is a critical pair
$\af,\xf$ such that $f(\xf,\af,3)$ is an inflexion point. Solving the pair of equations
$f(\xf,\af,4)=f(\xf,\af,5)=0$ gives $\af\approx 2.903$, $\xf\approx 1.344$.
In general, the function $f(x,a,n)$ will be monotonic decreasing for $a$ `small enough'.
From Alzer and Berg~\cite{Alzer:2002}, $a=1$ is always such a value.
The change to non-monotonic behaviour will occur at a value $\af(n)$ when $f(x,a,n)$ has
an inflexion point at $\xf(n)$. The values of $\af(n)$ and $\xf(n)$ can be determined
from the equations
\[f(\xf,\af,n+1)=f(\xf,\af,n+2)=0\ .\]

We next consider how $\af(n)$ depends upon $n$.
We have seen already that
$\af(2)\approx 3.138$, and $\af(3)\approx 2.903$.
For all $n>2$, $\af(n+1)<\af(n)$. This is so because if $f(x,\af(n),n)$ has an inflexion
point, then $f(x,\af(n),n+1)$ must be non-monotonic. Therefore, if one thinks of
$a$ as increasing from $a=1$, for which value all functions are monotonic, then $f(x,a,n+1)$
must have become non-monotonic at a lower value of $a$ than that for which $f(x,a,n)$
becomes non-monotonic.
Thus the sequence $\{\af(n)\}$ must be monotonically decreasing and it is bounded below by
$1$ from the results in \cite{Alzer:2002}.
Therefore the limit $\ac=\lim_{n\to\infty}\af(n)$
exists and is the critical value of $a$ such that $H(x;a)$ is
completely monotonic for $a < \ac$.

\section{Numerical method}
In view of the above observations, we need to compute $\af(n)$ for a selection of
values of $n$ and then extrapolate to compute the limit $n\to\infty$.
In order to obtain a value of $\ac$ that is accurate enough to be submitted
to tools such as Plouffe's inverter~\cite{Plouffe}, we need to compute derivatives up to the order of $10^5$.
Hence, an efficient method must be found to evaluate $f(x,a,n)$ for large values of $n$. Symbolic differentiation becomes impossible after at most $n=50$, depending upon the
memory available on the computer used, implying the need for a numerical scheme.

By introducing the notation
\[   h = \left(1+\frac{a}{x}\right)^{x}\ ,  \]
we can write $H(x;a)=\exp(a)-h$, and obtain
\begin{eqnarray}
  f(x,a,0) &=& H(x;a)= e^a - h \ , \\
  f(x,a,1) &=& \frac {\partial h}{\partial x} = \left(
 \ln  \frac {x+a}{x} -\frac {a}{x+a} \right)h \ .
\end{eqnarray}
Leibnitz's rule now allows us to compute recursively.
\begin{eqnarray}
f(x,a,n) &=&  (-1)^{n-1} \frac{\partial^{n-1}}{\partial x^{n-1}} f(x,a,1)  \nonumber \\
  &=& (-1)^{n-1}\sum_{k=0}^{n-1} \binom{n-1}{k} \frac{\partial^{k}h}{\partial x^{k}}
  \frac{\partial^{n-1-k}}{\partial x^{n-1-k}}\left( \ln  \frac {x+a}{x} -\frac {a}{x+a} \right) \label{eq:recurr} .
\end{eqnarray}
The second term in the sum can be computed explicitly by introducing the following function.
\[ g(x,a,n) =
      \begin{cases}
       -\ln  \frac {x+a}{x} +\frac {a}{x+a} \ ,& n=0\ ,\\
     -\frac{(n-1)!}{x^n}+\frac{(n-1)!}{(x+a)^n}
+\frac{n!a}{(x+a)^{n+1}}\ , & n\ge 1. \end{cases}
\]
Thus the computational form of (\ref{eq:recurr}) becomes
\begin{equation}
f(x,a,n) = \sum_{k=0}^{n-1} \binom{n-1}{k} f(x,a,k) g(x,a,n-1-k)
  \label{eq:comp}\ .
\end{equation}
Although it is possible to compute with this formula symbolically, expression swell
prevents this approach from being useful for anything other than for error checking.
If we assign numerical
values to $x$ and $a$ and compute $f(x,a,0)$ numerically, then (\ref{eq:comp}) can evaluate derivatives numerically to orders $10^6$ or more.

Once the functions $f(x,a,n)$ can be computed efficiently, the equations
$f(\xf,\af,n+1)=0$, $f(\xf,\af,n+2)=0$ can be solved to find the transition point in
$f(x,a,n)$. For small values of $n$, any method, for example Maple's \texttt{fsolve},
can be used. For large $n$, extended precision must be used, owing to loss of precision
through the accumulaton of cancellation and rounding errors.
For $n$ greater than $10^4$ the precision loss amounts to about 15 decimal digits.
The equations were solved using bivariate Newton iteration. The required derivatives
with respect to $a$ can be computed using a method analogous to that above.
In addition to cancelation errors, there is another reason for requiring extended precision. Because $(\xf,\af)$ is a repeated root of $f(\xf,\af,n+1)$, Newton iteration converges slowly and a nearly singular matrix must be inverted. Therefore 60 decimal digits were used in the computations. For smaller $n$, Maple was used, but became too slow for
larger $n$, and Bailey's \texttt{ARPREC}~\cite{Bailey} and Gnu GMP were used.

\section{Results}

Table \ref{tab:one} displays critical values computed between $n=1000$ and $n=100000$.
Each entry gives the smallest value of $a$ such that $f(\xf,a,n+1)<0$.
That is, for each entry, $f(\xf,\af-10^{-14},n+1)>0$.

\begin{table}
\begin{center}
\begin{tabular}{r|c|l}
n           & $\af(n)$ & $\xf(n)$ \\
\hline
       1000& 2.30183958971854&  436.380167908055 \\
       2000& 2.30074838075010&  872.743540008136 \\
       4000& 2.30020250313093&  1745.47034071250 \\
       5000& 2.30009330574014&  2181.83374672043 \\
      10000& 2.29987488908353&  4363.65078807689 \\
      12500& 2.29983120225165&  5454.55931101894 \\
      16000& 2.29979297531646&  6981.83124393025 \\
      20000& 2.29976566981568&  8727.28488210931 \\
      40000& 2.29971105744643&  17454.5530758350 \\
      50000& 2.29970013475374&  21818.1871732641 \\
     100000& 2.29967828914950&  43636.3576615316 \\
\end{tabular}
\end{center}
\caption{\label{tab:one} Table of critical values.}
\end{table}

\begin{figure}
\begin{center}
    \scalebox{0.6}{\includegraphics{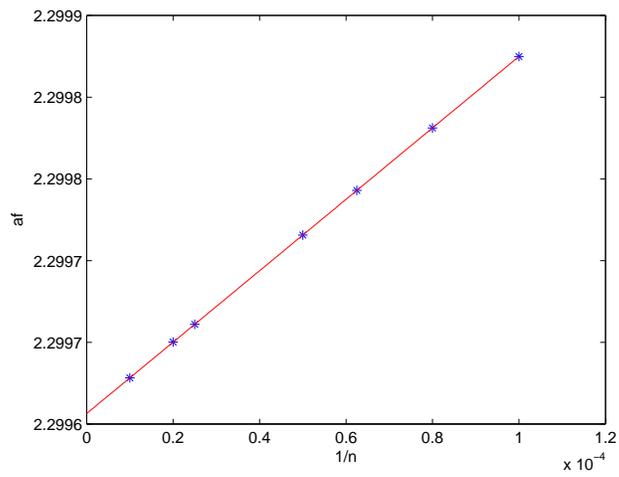}}
    \end{center}
    \caption{\label{fig:extr} A plot of $\af$ against $n^{-1}$.}
\end{figure}

To extrapolate to $n\to\infty$, we examine $\af$ as a function of
$n^{-1}$. Figure \ref{fig:extr} shows a plot of $\af(n)$ against $n^{-1}$.
We conjecture that $\af$ obeys a relation
\begin{equation}
      \af = \ac + \frac{a_1}{n} + \frac{a_2}{n^2} + O(n^{-3}).
\end{equation}
By fitting a quadratic in $n^{-1}$ using the data points $n=40000,\, 50000,\, 100000$,
we obtain $\ac \approx 2.29965644325$, where the error can be expected to be $O(10^{-12})$.
This value was given to Plouffe's inverter~\cite{Plouffe}, the inverse symbolic calculator~\cite{InvSym}, and Maple's \texttt{identify} command. No symbolic quantity matched all
digits, and the closest symbolic quantities offered no immediate inspiration for analytic investigations.

\section{Closing remarks}
The simple solution to problem (\ref{eq:known}) has led to
several attempts to find an equally simple solution for (\ref{eq:def}).
The numerical evidence presented suggests that such attempts will continue to fail.
The numerical data obtained here does, however, suggest a number of simple properties
which might help in an analytical
solution to the problem. The dependence of $\xf(n)$ on $n$ suggests a strongly linear correlation between $\xf$ and $n$. A similar correlation between $\ac(n)$ and $n^{-1}$
has been noted above. If these numerical observations can be explained analytically, then
a full solution might follow.

\end{document}